\documentclass[11pt]{article}
\usepackage[applemac]{inputenc}

\usepackage{epsf}
\usepackage{color}
\usepackage{epsfig}

\usepackage{graphicx}
\usepackage{amsmath}
\usepackage{latexsym}
\usepackage{amssymb}

\usepackage{alltt}

\usepackage{a4}

\newtheorem{theorem}{Theorem}[section]

\numberwithin{equation}{section}

\newcommand{\ad}{{\rm ad}}

\newcommand{\e}{{\rm e}}
\newcommand{\R}{R}

\begin{document}

\title{On the structure and convergence of the symmetric Zassenhaus formula
}

\author{Ana Arnal\thanks{Email: \texttt{ana.arnal@uji.es}}
 \and
 Fernando Casas\thanks{Corresponding author. Email: \texttt{Fernando.Casas@uji.es} }\and Cristina Chiralt\thanks{Email: \texttt{chiralt@uji.es}}
  }

\date{}
\maketitle

\begin{abstract}

We propose and analyze a symmetric version of the Zassenhaus formula for disentangling the exponential of two non-commuting operators. A recursive procedure for generating the
expansion up to any order is presented  which also allows one to get an enlarged domain of convergence when it is formulated for matrices. It is shown that
the approximations obtained by truncating the infinite expansion considerably improve those arising from the standard Zassenhaus formula.

\vspace*{1cm}

\begin{center}
Institut de Matem\`atiques i Aplicacions de Castell\'o (IMAC) and  De\-par\-ta\-ment de
Ma\-te\-m\`a\-ti\-ques, Universitat Jaume I,
  E-12071 Cas\-te\-ll\'on, Spain.
\end{center}

\end{abstract}

\section{Introduction}

It is a well known fact that, given two arbitrary $n \times n$ matrices $A$ and $B$, the following two statements are equivalent:
\begin{itemize}
  \item[(a)] $A$ and $B$ commute.
  \item[(b)] $\exp(\lambda A + \lambda B) = \exp(\lambda A) \, \exp(\lambda B)$ for all values of the scalar $\lambda$.
\end{itemize}
If, on the other hand, $A B \ne B A$ and the matrix $\e^{\lambda A} \, \e^{\lambda B}$ is the exponential of some matrix $D$, then $D$ can be written as
$D = \lambda(A + B) + C$,  
where the additional term $C$ is due to the non-commutativity of $A$ and $B$. This, of course, is related with
 the celebrated Baker--Campbell--Hausdorff formula \cite{bonfiglioli12tin}.

 One could also consider the dual problem, namely, to get matrices $C_1, C_2, \ldots$ such that $\e^{\lambda(A+B)} = \e^{\lambda A} \, \e^{\lambda B} \, \e^{C_1(\lambda)} \,
 \e^{C_2(\lambda)} \cdots$, so that the exponential $\e^{\lambda(A+B)}$ can be `disentangled'. The Zassenhaus
 formula establishes this fact that provides explicit expressions for $C_1, C_2, \ldots$ in terms of $A$ and $B$ \cite{magnus54ote}. It turns out that the commutator
 \[
   [A,B] = AB - BA
 \]
 plays a fundamental role when analyzing both issues.

 To establish the Zassenhaus formula in a general setting, we consider two  non commuting indeterminate variables $X$, $Y$ and the unital associative algebra
 of formal power series in $X$, $Y$ over the field $\mathbb{K}$, denoted by $\mathbb{K}\langle X, Y \rangle$.  Here $\mathbb{K}$ can be, say, $\mathbb{Q}$, $\mathbb{R}$
 or $\mathbb{C}$. With the operation
 $X, Y \longmapsto [X, Y]$ it is possible to introduce a commutator algebra $[\mathbb{K}\langle X, Y \rangle ]$. Then, in $[\mathbb{K}\langle X, Y \rangle ]$ the set
 of all Lie polynomials in $X$ and $Y$ (i.e., all possible expressions obtained from $X$, $Y$ by addition, multiplication by elements in $\mathbb{K}$ and the
 commutator) forms a subalgebra $\mathcal{L}(X,Y)$, which turns out to be a free Lie algebra with generators $X$, $Y$  \cite{postnikov94lga,reutenauer93fla}.
 It is clear then that every Lie polynomial in  $\mathcal{L}(X,Y)$ is also a polynomial in $\mathbb{K}\langle X, Y \rangle$ (one only has to remove all commutators
 by the rule $[a,b] = ab - ba$). 
 
 With this notation, the Zassenhaus formula establishes that the exponential $\e^{X+Y}$ can be uniquely decomposed as
\begin{equation}  \label{zass.1}
      \e^{X + Y} = \e^X \, \e^Y \, \prod_{n=2}^{\infty} \e^{C_n(X,Y)} = \e^X \, \e^Y \, \e^{C_2(X,Y)} \,
      \e^{C_3(X,Y)} \, \cdots  \, \e^{C_k(X,Y)} \, \cdots,
\end{equation}
where $C_k(X,Y)$ is a homogeneous Lie polynomial in $X$ and $Y$ of degree $k$
\cite{magnus54ote,suzuki77otc,weyrauch09ctb,wilcox67eoa,witschel75ooe}. The first terms read explicitly
 \[
\aligned
  &  C_2(X,Y) = -\frac{1}{2} [X,Y] \\
  &  C_3(X,Y) = \frac{1}{3} [Y,[X,Y]] + \frac{1}{6} [X,[X,Y]] \\
  &  C_4 = -\frac{1}{24} [X,[X,[X,Y]]] - \frac{1}{8} [Y,[X,[X,Y]]] - \frac{1}{8} [Y,[Y,[X,Y]]].
\endaligned
\]
There are several procedures to get the terms in the expansion (\ref{zass.1}) \cite{suzuki77otc,weyrauch09ctb,witschel75ooe}. In particular,
a recursive algorithm has been designed in  \cite{casas12eco} for obtaining
$C_n$ up to a prescribed value of $n$ directly in terms of the minimum number of independent commutators involving $n$ operators $X$ and
$Y$. The algorithm, in addition, can be easily implemented in a symbolic
algebra system without any special requirement, beyond the linearity property of the commutator.

Given the ubiquity of exponentials of sums of operators in physics and mathematics, it is natural to find applications of  Zassenhaus formula  in many
different fields: periodically driven quantum systems \cite{goldman14pdq}, quantum nonlinear optics \cite{quesada14eot}, $q$-analysis in quantum groups \cite{quesne04dqe},
 the Schr\"odinger equation in the semiclassical regime \cite{bader14eaf}, hypoelliptic differential equations \cite{hormander67hso}
 and splitting methods in numerical analysis \cite{geiser11osm}, just to quote a few.
 Of course, in these applications one truncates the infinite product (\ref{zass.1}) at some $n$ and then takes the approximation
 \[
     \e^{X + Y} \approx \e^X \, \e^Y e^{C_2(X,Y)} \, \e^{C_3(X,Y)} \, \cdots  \, \e^{C_n(X,Y)}.
\]

 Whereas the factorization (\ref{zass.1}) is well defined in the above setting, i.e. in $\mathbb{K}\langle X, Y \rangle$ (or more specifically, in the 
 subalgebra $\mathcal{L}(X,Y)$ of $[\mathbb{K}\langle X, Y \rangle ]$), it has only a finite radius of convergence when
$X$ and $Y$ are elements of a Banach Lie algebra (in particular, when they are
 $n \times n$ real or complex matrices). Specifically,
 \begin{equation}   \label{conv.1}
   \lim_{n \rightarrow \infty} \e^{X} \, \e^{Y} \,
   \e^{C_2} \cdots \e^{C_n} = \e^{X+Y}
\end{equation}
only in a certain subset of the plane $(\|X\|, \|Y\|)$ \cite{suzuki77otc,bayen79otc}.
Thus, in \cite{casas12eco} it is shown that the convergence domain contains the region
$\|X\|+\|Y\| < 1.054$, and extends to the points $(\|X\|,0)$ and $(0,\|Y\|)$
with arbitrarily large values of $\|X\|$ or $\|Y\|$.

Instead of (\ref{zass.1}), it is also possible to consider the ``left-oriented" Zassenhaus formula
\begin{equation}   \label{zass.1.1}
    \e^{X + Y} = \cdots \, \e^{\hat{C}_4(X,Y)} \, \e^{\hat{C}_3(X,Y)} \, \e^{\hat{C}_2(X,Y)} \,
        \e^Y \, \e^X,
\end{equation}
but the respective exponents are closely related: $\hat{C}_i(X,Y) = (-1)^{i+1} C_i(X,Y)$,  for $i \ge 2$. A third possibility would consist in taking
left-right palindromic factorizations, either of the form
\begin{equation}  \label{zass.1.2}
    \e^{X + Y} = \e^{\frac{X}{2} } \, \e^{\frac{Y}{2}  } \,
   \e^{C_2} \, \e^{ C_3} \,  \cdots  \e^{ C_n}  \cdots\,  \e^{ C_n} \cdots \, \e^{ C_3} \,
   \e^{ C_2} \, \e^{\frac{ Y}{2} } \, \e^{\frac{X}{2} }
\end{equation}
or
\begin{equation}  \label{zass.1.2b}
    \e^{X + Y} = \cdots \e^{ D_n} \, \cdots \e^{ D_2} \, \e^{\frac{Y}{2}} \, \e^X \, \e^{\frac{Y}{2}} \, \e^{ D_2} \cdots \e^{ D_n} \cdots.
\end{equation}
 As a matter of fact, left-right symmetric compositions are usually preferable when integrating numerically systems of differential equations with
splitting methods \cite{mclachlan02sm}. In this context,  the leapfrog/St\"ormer--Verlet/Strang splitting is a paradigmatic example.
These methods preserve the time-symmetry of the continuous system and the asymptotic
expansion of their local error only contains even powers of the step size parameter, thus conferring favorable properties to the numerical approximations \cite{blanes16aci,hairer06gni}.

Motivated by the recent use of \cite{casas12eco} in different contexts (i.e., periodically driven quantum systems \cite{goldman14pdq}, quantum predictive filtering \cite{grimsmo15qpf}, motion of a quantum particle in a
magnetic monopole field \cite{soloviev16dmq}, quantum nonlinear optics \cite{quesada14eot}) and the aforementioned advantages of left-right symmetric compositions, in this paper we 
explore the possibility of constructing in a systematic way the different exponents in the factorization (\ref{zass.1.2}) and test its main features on several examples. More
specifically, we have the following goals in mind: first, we wish to establish its validity; second, to enlarge the convergence domain of the usual Zassenhaus formula and also to speed up
its rate of convergence, and third, to discern whether it is possible to get a more accurate approximation when the factorization is truncated at some given term $n$.
In fact, approximations of this kind have already been considered by \cite{bader14eaf,bader16emf} 
in the context of designing new integration methods for the
Schr\"odinger equation in the semiclassical regime, whereas in \cite{grimsmo15qpf} the composition (\ref{zass.1.2b}) is used for Hermitian operators. 

Here we propose a
systematic procedure that allows us to get the terms $C_k$ in the factorization (\ref{zass.1.2}) in a more efficient way than by considering the Baker--Campbell--Hausdorff formula
and also to 
obtain better estimates for the convergence of the expansion. Finally, to fully appreciate the advantages of using the symmetric version of the Zassenhaus formula 
(\ref{zass.1.2}) instead of the usual one (\ref{zass.1}), we apply it to a pair of physical examples involving finite-dimensional Lie subalgebras of operators, as well as real matrices.

\section{A symmetric version of the Zassenhaus formula}

In this section we constructively prove the following theorem, which constitutes a generalization of the factorization (\ref{zass.1}).
\begin{theorem}   \label{theo.1}
  (Symmetric Zassenhaus Formula). Under the same hypotheses than the usual Zassenhaus formula (\ref{zass.1}), $\e^{X + Y}$ can also be uniquely decomposed  as
  \begin{eqnarray}  \label{zassym}
      \e^{X + Y} & = &  \e^{\frac{X}{2} } \, \e^{\frac{Y}{2}  } \left( \prod_{n=2}^{\infty} \e^{C_n(X,Y)} \right) \,
       \left( \prod_{n=\infty}^{2} \e^{C_n(X,Y)} \right)  \, \e^{\frac{ Y}{2} } \, \e^{\frac{X}{2} }  \nonumber  \\
    & \equiv & \e^{\frac{X}{2} } \, \e^{\frac{Y}{2}  }  \e^{C_2} \, \e^{ C_3} \,  \cdots  \e^{ C_k}  \cdots\, \e^{ C_k}\, \cdots \e^{ C_3} \,
   \e^{ C_2} \, \e^{\frac{ Y}{2} } \, \e^{\frac{X}{2} } ,
  \end{eqnarray}
  where $C_k(X,Y)$ is again a homogeneous Lie polynomial in $X$ and $Y$ of degree $k$ (in principle different from the corresponding Lie polynomials appearing
  in (\ref{zass.1})). Moreover, in (\ref{zassym}) all terms
  $C_{2k}$ vanish identically, so that
    \begin{equation}  \label{zassym.red}
      \e^{X + Y} =  \e^{\frac{X}{2} } \, \e^{\frac{Y}{2}  } \,
   \e^{C_3} \, \e^{ C_5} \,  \cdots  \e^{ C_{2k+1}}  \cdots\, \e^{ C_{2k+1}}\, \cdots \e^{ C_5} \,
   \e^{ C_3} \, \e^{\frac{ Y}{2} } \, \e^{\frac{X}{2} }.
  \end{equation}
\end{theorem}

\

We first show that $C_{2k} \equiv 0$ for all $k$, i.e., only Lie polynomials of odd degree appear in (\ref{zassym}).
To see this, let us introduce a parameter $\lambda > 0$ multiplying each variable
$X$ and $Y$ and denote
\begin{eqnarray}  \label{zass.3.1}
   \Psi(\lambda) & \equiv & \e^{\lambda (X+Y)}  \\
    & = & \e^{\frac{\lambda}{2} X} \, \e^{\frac{\lambda}{2}  Y} \,
   \e^{\lambda^2 C_2} \, \e^{\lambda^3 C_3} \,  \cdots  \e^{\lambda^k C_k}  \cdots\, \e^{\lambda^k C_k}\, \cdots \e^{\lambda^3 C_3} \,
   \e^{\lambda^2 C_2} \, \e^{\frac{\lambda}{2}  Y} \, \e^{\frac{\lambda}{2} X}. \nonumber
\end{eqnarray}
On the one hand it is true that $\Psi^{-1}(-\lambda) = \big(\e^{-\lambda (X+Y)}\big)^{-1} = \Psi(\lambda)$. On the other hand,
\[
 \Psi^{-1}(-\lambda) =  \e^{\frac{\lambda}{2} X} \, \e^{\frac{\lambda}{2}  Y} \,
   \e^{-\lambda^2 C_2} \, \e^{\lambda^3 C_3} \,    \e^{-\lambda^4 C_4}  \cdots\, \e^{-\lambda^4 C_4}\,  \e^{\lambda^3 C_3} \,
   \e^{-\lambda^2 C_2} \, \e^{\frac{\lambda}{2}  Y} \, \e^{\frac{\lambda}{2} X},
\]
so that by comparing with the second line of (\ref{zass.3.1}), it follows that all terms containing an even number of operators $X$, $Y$ in (\ref{zassym}) vanish
identically.

\

The terms $C_k$ in (\ref{zass.3.1}) can be obtained, of course, in a number of ways. For instance, by taking the $n$-derivative of both terms in
(\ref{zass.3.1}) with respect to $\lambda$ at $\lambda=0$ and equating both of them lead to the expression of $C_n$, $n=2,3,\ldots$. Alternatively, one has
\begin{eqnarray}   \label{suzu.1}
      C_n & = & \frac{1}{2 n!} \Big( \frac{d^n }{d \lambda^n} \big(  \e^{-\lambda^{n-1} C_{n-1}}
      \cdots \e^{-\lambda^2 C_2} \e^{-\frac{\lambda}{2} Y} \e^{-\frac{\lambda}{2} X}
         \e^{\lambda(X+Y)}   \nonumber \\
      &  &   \quad  \e^{-\frac{\lambda}{2} X}   \e^{-\frac{\lambda}{2} Y} \e^{-\lambda^2 C_2} \cdots   \e^{-\lambda^{n-1} C_{n-1}} \big) \Big)_{\lambda=0}.
 \end{eqnarray}
Working out this recursion is not an easy task, however, and furthermore it does not provide $C_n$ directly
in terms of elements in $\mathcal{L}(X,Y)$, i.e., in terms of commutators. Another possibility consists in applying the
Baker--Campbell--Hausdorff formula, as in \cite{bader14eaf}, but of course one needs to construct this formula in advance and the procedure is not optimal
with respect to computer time and memory resources.

Next, we present yet another recursive algorithm for generating $C_n$ directly in terms of commutators.
The procedure can be better understood by considering the general
composition (\ref{zassym}). This, in addition, can be used as an additional check for the algorithm, in the sense that all the exponents $C_{2k}$
obtained with it have to vanish.

First we introduce the  products
\begin{equation}  \label{foz.1}
\aligned
  & R_1(\lambda) =  \e^{-\frac{\lambda}{2} Y} \, \e^{-\frac{\lambda}{2} X} \, \e^{\lambda(X+Y)},  \\
  & R_n(\lambda)  =  \e^{-\lambda^{n} C_{n}} \, \e^{-\lambda^{n-1} C_{n-1}} \,
      \cdots \, \e^{-\lambda^2 C_2} \, R_1(\lambda)  = \e^{-\lambda^{n} C_{n}} \, R_{n-1}(\lambda),  \qquad n \ge 2.
\endaligned
\end{equation}
It is clear from (\ref{zass.3.1}) that
\begin{equation}   \label{foz.2}
   R_n(\lambda) =  \e^{\lambda^{n+1} C_{n+1}}  \, \e^{\lambda^{n+2} C_{n+2}} \cdots \e^{\lambda^{n+2} C_{n+2}} \cdots
   \e^{\lambda^2 C_2} \ \e^{\frac{\lambda}{2} Y} \,   \e^{\frac{\lambda}{2} X}
\end{equation}
for all $n \ge 1$.
In this expression, as in (\ref{zassym}), the first dots indicate that the index $n$ first increases (up to infinity) and then decreases again.
Finally, we introduce
  \begin{equation}
\label{eq:F_nDef}
    F_n(\lambda) \equiv  \left(
\frac{d}{d \lambda} R_n(\lambda)
\right) R_n(\lambda)^{-1}, \qquad n \ge 1.
  \end{equation}
When $n=1$, and taking into account the expression of $R_1(\lambda)$ given in (\ref{foz.1}), we get
\begin{equation}  \label{eq.F1}
 F_1(\lambda)
=\left(
\frac{d}{d \lambda} \R_1(\lambda)
\right) \R_1(\lambda)^{-1}=
-\frac{Y}{2}+\e^{-\frac{\lambda}{2} \ad_{Y}}  \left( \e^{-\frac{\lambda}{2} \ad_{X}}Y+\frac{X}{2} \right),
\end{equation}
where
\[
  \e^{\ad_A} B =  \e^{A} B \e^{-A}   = \sum_{n\ge 0} \frac{1}{n!} \ad_A^n B
\]
and the ``ad" operator is defined as
\[
\ad_A B = [A,B], \qquad  \ad_A^j B = [A, \ad_A^{j-1} B], \qquad \ad_A^0 B = B.
\]
Working out expression (\ref{eq.F1}) one arrives at the power series
\begin{equation}
  \label{eq:F_1}
 F_1(\lambda) = \sum_{\ell=0}^{\infty}  f_{1,\ell} \, \lambda^{\ell},
\end{equation}
with
\begin{eqnarray}   \label{eq:F_1.1}
 f_{1,0} & = & \frac{1}{2} (X + Y), \\
 f_{1,\ell} & = & \frac{(-1)^{\ell}}{2^{\ell}} \left( \frac{1}{2 \ell!} \ad_{Y}^{\ell}X+ \sum_{j=0}^{\ell}\,\frac{1}{j!(\ell-j)!}\; \ad_{Y}^{\ell-j}\; \ad_{X}^{j}Y \right),
    \qquad \ell \ge 1. \nonumber
\end{eqnarray}
A similar expansion  can be obtained for $F_n(\lambda)$, $n \ge 2$, by considering  the expression of $R_n(\lambda)$ given in
(\ref{foz.1}) and the relation (\ref{eq:F_nDef}). More specifically,
\begin{equation}  \label{eq:F_n.n}
\aligned
  & F_n(\lambda) =  -n \, C_n \, \lambda^{n-1}
+ \e^{-\lambda^n C_n} \, \left( \frac{d }{d \lambda} R_{n-1}(\lambda) \right) R_{n-1}(\lambda)^{-1} \,  \e^{\lambda^n C_n}  \\
   &  \;\; =    -n \, C_n \, \lambda^{n-1}
+ \e^{-\lambda^n C_n} \, F_{n-1}(\lambda) \, \e^{\lambda^n C_n}
=  -n \, C_n \, \lambda^{n-1}+ \e^{-\lambda^n \ad_{C_n}} F_{n-1}(\lambda)    \\
&   \;\; =  \e^{-\lambda^n \ad_{C_n}} (F_{n-1}(\lambda) -n \, C_n \, \lambda^{n-1}),
\endaligned
\end{equation}
so that
   \begin{equation}   \label{rec.2.3}
     F_n(\lambda) = \sum_{\ell=0}^{\infty}  f_{n,\ell} \, \lambda^{\ell}, \qquad \mbox{ with } \qquad
        f_{n,\ell} = \sum_{j=0}^{[\ell/n]}  \frac{(-1)^j}{j!} \ad_{C_n}^{j} \widetilde{f}_{n-1,\ell-nj},
 \end{equation}
where $[\ell/n]$ denotes the integer part of $\ell/n$ and
\[
   \widetilde{f}_{n-1,\ell} = \left\{  \begin{array}{lc}
   		f_{n-1,\ell} \  & \ell \ne n-1 \\
		f_{n-1,\ell} - n C_n, \ &  \ell = n-1.
	\end{array} \right.
\]		
If we take instead the expression (\ref{foz.2}) for $R_1$ and evaluate again $F_1$ according with (\ref{eq:F_nDef}) we arrive at
\begin{eqnarray}  \label{f1.2}
   F_1(\lambda) & = &  2 \, C_{2} \, \lambda + \sum_{j=3}^{\infty} j \, \lambda^{j-1} \,
  \e^{\lambda^{2} \ad_{C_{2}}} \cdots \e^{\lambda^{j-1} \ad_{C_{j-1}}} C_j  \nonumber \\
   &  & +\sum_{j=2}^{\infty} j \, \lambda^{j-1} \, \left( \prod_{k=2}^{\infty} \e^{\lambda^k \ad_{C_k}} \right)
     \left( \prod_{k=\infty}^{j+1} \e^{\lambda^k \ad_{C_k}} \right)  C_j  \\
   &  & +  \left( \prod_{k=2}^{\infty} \e^{\lambda^k \ad_{C_k}} \right) \left( \prod_{k=\infty}^{2} \e^{\lambda^k \ad_{C_k}} \right)
   G_1(\lambda),  \nonumber
\end{eqnarray}
where
\[
   G_{1}(\lambda) =\dfrac 12\left( Y+ \e^{\frac \lambda 2 \ad_{Y}}X\right).
\]
Explicitly,
\begin{equation}
  \label{eq:g_1}
 G_1(\lambda) = \sum_{\ell=1}^{\infty} g_{1,\ell} \,  \lambda^{\ell},
\end{equation}
with
\begin{equation}   \label{eq:G_1.1}
 g_{1,0}  =  \frac{1}{2} (X + Y), \qquad\quad
 g_{1,\ell}  =  \frac{1}{\ell!\;2^{\ell+1}} \;\ad_{Y}^{\ell}X.
    \qquad \ell \ge 1.
\end{equation}
Analogously, for $n \ge 2$ one has
\begin{eqnarray}  \label{fn.2}
   F_n(\lambda) & = &  (n+1) \, C_{n+1} \, \lambda^n + \sum_{j=n+2}^{\infty} j \, \lambda^{j-1} \,
  \e^{\lambda^{n+1} \ad_{C_{n+1}}} \cdots \e^{\lambda^{j-1} \ad_{C_{j-1}}} C_j  \nonumber \\
   &  &  +\sum_{j=n+1}^{\infty} j \, \lambda^{j-1} \, \left( \prod_{k=n+1}^{\infty} \e^{\lambda^k \ad_{C_k}} \right)
     \left( \prod_{k=\infty}^{j+1} \e^{\lambda^k \ad_{C_k}} \right)  C_j  \\
   &  & +  \left( \prod_{k=n+1}^{\infty} \e^{\lambda^k \ad_{C_k}} \right) \left( \prod_{k=\infty}^{n+1} \e^{\lambda^k \ad_{C_k}} \right)
   G_n(\lambda),  \nonumber
\end{eqnarray}
in terms of
\[
  G_{n}(\lambda) = n \lambda^{n-1} C_n+ \e^{\lambda^{n} \ad_{C_{n}}} G_{n-1}(\lambda).
\]
In this case,
\begin{equation}   \label{rec.G_n}
     G_n(\lambda) = \sum_{\ell=0}^{\infty}  g_{n,\ell} \, \lambda^{\ell}, \qquad \mbox{ with } \qquad
        g_{n,\ell} = \sum_{j=0}^{[\ell/n]}  \frac{1}{j!} \ad_{C_n}^{j} \widetilde{g}_{n-1,\ell-nj},
 \end{equation}
and
\[
   \widetilde{g}_{n-1,\ell} = \left\{  \begin{array}{lc}
   		g_{n-1,\ell} \  & \ell \ne n-1 \\
		g_{n-1,\ell} +n C_n, \ &  \ell = n-1.
	\end{array} \right.
\]	
Now, comparing (\ref{eq:F_1}) with (\ref{f1.2}), it is clear that the term independent of $\lambda$ in both expressions is exactly the same
($f_{1,0} = g_{1,0}$), whereas for the term in $\lambda$ we get
\[
  f_{1,1} = 2 C_2 + 2 C_2 + g_{1,1},
\]
so that
\[
  C_2 = \frac{1}{4} (f_{1,1} - g_{1,1}) = \frac{1}{4} \left(   - \frac{1}{4} \ad_Y X - \frac{1}{2} \ad_X Y  - \frac{1}{4} \ad_Y X \right) = 0.
\]
The exponent $C_3$ can be obtained by comparing analogously terms in $\lambda^2$ in the expressions (\ref{rec.2.3}) and (\ref{fn.2})
when $n=2$: $f_{2,2} = 3 C_3 + 3 C_3 + g_{2,2}$. Proceeding by induction, we get in general for the term in $\lambda^{k-1}$
\[
  f_{k-1,k-1} = 2 k C_k + g_{k-1,k-1}
\]
whence
\begin{equation}
\label{eq:Cn}
C_{k} = \frac{1}{2k} \, \left(f_{k-1,k-1}-g_{k-1,k-1}\right), \qquad k \ge 2.
\end{equation}

From the previous discussion, the algorithm for generating the exponents $C_{k}$, $k=3,5,7,\ldots$ in (\ref{zassym.red}) can be
formulated as follows.
\begin{equation}   \label{alg.1}
\begin{array}{l}
  \mbox{Define} \; f_{1,k} \; \mbox{by eq. (\ref{eq:F_1.1})}, \; k \ge 0  \\
  \mbox{Define} \; g_{1,k} \; \mbox{by eq. (\ref{eq:G_1.1})}, \; k \ge 0  \\
   C_3 = \frac{1}{6} (f_{1,2} - g_{1,2}) \\
   \mbox{For } k=3,5,\ldots  \\
   \quad    f_{k,0} = f_{1,0} \\
   \quad    g_{k,0} = g_{1,0} \\
   \quad     f_{k, \ell} = \sum_{j=0}^{[\ell/k]} \frac{(-1)^j}{j!} \ad_{C_k}^j f_{k-2,\ell - k j} \\
   \quad     f_{k,k-1} = - k C_k + f_{k-2,k-1} \\
   \quad     g_{k, \ell} = \sum_{j=0}^{[\ell/k]} \frac{1}{j!} \ad_{C_k}^j g_{k-2,\ell - k j} \\
   \quad     g_{k,k-1} =  k C_k + g_{k-2,k-1} \\
   \quad     C_{k+2} = \frac{1}{2(k+2)} (f_{k,k+1} - g_{k,k+1})
\end{array}
\end{equation}
Here we must impose that $C_{2k}=0$.  This procedure can be easily implemented in a symbolic algebra system to render explicit
expressions for the terms $C_{2k+1}$ up to any value of $k$, although, in general, there is no guarantee \textit{a priori}
that all commutators are independent, in contrast with the
 standard Zassenhaus formula \cite{casas12eco}. Of course, algorithm (\ref{alg.1}) can be combined with the technique presented in \cite{casas09aea}
 for constructing a Hall basis in the free Lie algebra $\mathcal{L}(X,Y)$ so as to get the exponents  in that basis. In any case,
 up to $k=6$  it turns out that the  $C_{2k+1}$ produced as output by (\ref{alg.1}) are more compact than in the Hall basis, with a reduced
 computational cost. For illustration, the number of terms in $C_3, \ldots, C_{13}$ is: $2, 6, 18, 54,
 132, 630$ (algorithm (\ref{alg.1})) and $2,6,18,56,186,630$ (in the classical Hall basis). Notice in particular the remarkable reduction of terms in $C_{11}$.
 For illustration, the first terms read explicitly
\[
\aligned
 & C_3 =  \frac{1}{48}\, [X,[X,Y]]+ \frac{1}{24}\, [Y,[X,Y]] \\
 & C_5 =  \frac{1}{3840}\,[X,[X,[X,[X,Y]]]] + \frac{1}{960}\,[Y,[X,[X,[X,Y]]]]+\frac{1}{640}\,[Y,[Y,[X,[X,Y]]]]\\
&\ +\frac{1}{960}\,[Y,[Y,[Y,[X,Y]]]]-\frac{1}{960}\,[[X,Y],[X,[X,Y]]]-
\frac{1}{480}\,[[X,Y],[Y,[X,Y]]].
\endaligned
\]
A simple remark is worth noticing here: we can interchange the role of $X$ and $Y$ and Theorem \ref{theo.1} is still valid, i.e., it is true that
   \begin{equation}  \label{zassym.red2}
      \e^{X + Y} =  \e^{\frac{Y}{2} } \, \e^{\frac{X}{2}  } \,
   \e^{\tilde{C}_3} \, \e^{\tilde{C}_5} \,  \cdots  \e^{ \tilde{C}_{2k+1}}  \cdots\, \e^{ \tilde{C}_{2k+1}}\, \cdots \e^{ \tilde{C}_5} \,
   \e^{ \tilde{C}_3} \, \e^{\frac{ X}{2} } \, \e^{\frac{Y}{2} }.
  \end{equation}
Moreover, although the terms $\tilde{C}_i$ differ from those appearing in (\ref{zassym.red}), they can be generated by applying the same
recurrences (\ref{alg.1}) with the interchange $X \leftrightarrow Y$.

\section{Convergence of the symmetric Zassenhaus formula}

Algorithm (\ref{alg.1}) not only generates in an efficient way all the terms in the symmetric Zassenhaus formula (\ref{zassym}) but is also very useful for establishing its convergence when $X$ and $Y$ are elements of a Banach algebra $\mathcal{A}$, i.e.,
an algebra that is also a complete normed linear space with a sub-multiplicative norm,
\[
    \|X \, Y \| \le \|X\| \, \|Y\|.
\]
If this is the case, then $\| \ad_X Y = [X,Y]\| \le 2 \, \|X\| \, \|Y\|$ and, in general, $\| \ad_X^n Y\| \le 2^n \|X\|^n \, \|Y\|$.

To analyze the convergence, we proceed as in \cite{bayen79otc} and introduce, for $n \ge 2$, the truncated
left-right palindromic expansion
  \begin{equation}   \label{conv.2}
   \Psi_n(\lambda) \equiv \e^{\frac{\lambda}{2} X} \, \e^{\frac{\lambda}{2}  Y} \,
   \e^{\lambda^2 C_2} \, \e^{\lambda^3 C_3} \,  \cdots  \e^{\lambda^n C_n}  \, \e^{\lambda^n C_n}\, \cdots \e^{\lambda^3 C_3} \,
   \e^{\lambda^2 C_2} \, \e^{\frac{\lambda}{2}  Y} \, \e^{\frac{\lambda}{2} X}
\end{equation}
with $\lambda > 0$. For $n \ge 2$ and any $k>0$, one has
\begin{equation}   \label{conv.3}
  \Psi_{n+k}(\lambda) - \Psi_n(\lambda) = T_n(\lambda) \, U_{n+k}(\lambda) \, T_n^*(\lambda),
\end{equation}
where
\begin{eqnarray*}
   T_n(\lambda) & \equiv & \e^{\frac{\lambda}{2} X} \, \e^{\frac{\lambda}{2}  Y} \,
   \e^{\lambda^2 C_2} \, \e^{\lambda^3 C_3} \,  \cdots  \e^{\lambda^n C_n} \\
   U_{n+k}(\lambda) & \equiv &
    \e^{\lambda^{n+1} C_{n+1}} \cdots \e^{\lambda^{n+k} C_{n+k}} \,
   \e^{\lambda^{n+k} C_{n+k}} \cdots  \e^{\lambda^{n+1} C_{n+1}} - I
\end{eqnarray*}
and
\[
   T_n^*(\lambda)  = \e^{\lambda^n C_n} \cdots \e^{\lambda^3 C_3} \, \e^{\lambda^2 C_2} \,
     \e^{\frac{\lambda}{2}  Y} \,  \e^{\frac{\lambda}{2} X}.
\]
Now it is clear that
\[
   \|T_n(\lambda)\| \le \exp \left( \frac{1}{2} \lambda \|X\| +  \frac{1}{2} \lambda \|Y\| + \sum_{j=2}^n \lambda^j \|C_j\| \right)
\]
and an identical bound holds for  $\|T_n^*(\lambda)\|$.
On the other hand, from the estimate
\[
   \| \e^A \, \e^B - I \| \le \e^{\|A\| + \|B\|} -1,
\]
it follows by induction that
\[
  \|U_{n+k}(\lambda)\| \le \exp \left( 2 \sum_{j=n+1}^{n+k} \lambda^j \|C_j\| \right) -1.
\]
In consequence,
\[
\aligned
 &  \|\Psi_{n+k}(\lambda) - \Psi_n(\lambda)\| \le  \\
  &  \quad \exp \left(  \lambda (\|X\| +  \|Y\|) + 2 \sum_{j=2}^n \lambda^j \|C_j\| \right) \,
       \left( \exp \Big( 2 \sum_{j=n+1}^{n+k} \lambda^j \|C_j\| \Big) -1 \right).
\endaligned
\]
Suppose now that the series
\begin{equation}  \label{conv.exp}
    M(\lambda) \equiv  \lambda (\|X\| +  \|Y\|) + 2 \sum_{j=2}^{\infty} \lambda^j \|C_j\|
\end{equation}
has a certain radius of convergence $r_z$, say. Then, for $\lambda < r < r_z$ we have
\[
    2 \sum_{j=n+1}^{n+k} \lambda^j \|C_j\| \rightarrow 0 \quad \mbox{ as } \quad
n \rightarrow \infty
\]
and
\[
   \exp \left(  \lambda (\|X\| +  \|Y\|) + 2 \sum_{j=2}^n \lambda^j \|C_j\| \right) \le \e^{M(r)},
\]
so that
\[
   \|\Psi_{n+k}(\lambda) - \Psi_n(\lambda)\| \le \e^{M(r)} \,  \left( \exp \Big( 2 \sum_{j=n+1}^{n+k} r^j \|C_j\| \Big) -1 \right) \equiv V_n(r).
\]
It is then clear that $V_n(r) \rightarrow 0$ as $n \rightarrow \infty$, or equivalently
\[
            \|\Psi_{n+k}(\lambda) - \Psi_n(\lambda)\|  \rightarrow 0   \qquad \mbox{ as } \qquad n \rightarrow
            \infty,
\]
whence the sequence of entire functions $\Psi_n$ converges uniformly on any compact subset
of the ball $B(0,r_z)$, since $\mathcal{A}$ is complete. Let $\Psi$  denote this limit,
\[
    \lim_{n \rightarrow \infty} \Psi_n = \Psi.
\]
The function $\Psi(\lambda)$ is analytic in $B(0, r_z)$. Moreover, with the $C_k$ obtained with algorithm (\ref{alg.1}) (or alternatively, with expression
(\ref{suzu.1})),  one has
\[
   \left( \frac{d^k }{d \lambda^k} (  \Psi(\lambda ) \right)_{\lambda=0}
   =  \lim_{n \rightarrow \infty} \left( \frac{d^k }{d \lambda^k}
       (  \Psi_n(\lambda ) \right)_{\lambda=0}  = (X + Y)^k, \qquad \quad k \le n
\]
and thus $\Psi(\lambda)$ is indeed  $\e^{\lambda (X + Y)}$.   We have then shown that, for $\lambda \in B(0,r_z)$, the domain of convergence of the series
(\ref{conv.exp}), it is true that
\[
\lim_{n \rightarrow \infty} \left( \e^{\frac{\lambda}{2} X} \, \e^{\frac{\lambda}{2}  Y} \,
   \e^{\lambda^2 C_2} \, \e^{\lambda^3 C_3} \,  \cdots  \e^{\lambda^n C_n}  \, \e^{\lambda^n C_n}\, \cdots \e^{\lambda^3 C_3} \,
   \e^{\lambda^2 C_2} \, \e^{\frac{\lambda}{2}  Y} \, \e^{\frac{\lambda}{2} X} \right) = \e^{\lambda(X+Y)}
\]
and moreover, the convergence is uniform
on any compact subset of the ball $B(0, r_z)$. We are thus led to obtain an estimate for $r_z$, the radius of convergence of the power series (\ref{conv.exp}).

\

To this end, let us denote $\|X\| = x$, $\|Y\| = y$. From (\ref{eq:F_1.1}) we have $\|f_{1,0}\| \le \frac{1}{2} (x+y) \le (x+y)$ and
\begin{eqnarray}    \label{cota1:f1k}
 \|f_{1,\ell}\| & \leq &  \frac{1}{\ell!\;2^{\ell+1}}\,\|\ad_{Y}^{\ell}X\|+\frac{1}{2^\ell}\sum_{j=0}^{\ell}\,\frac{1}{j!(\ell-j)!}\; \|\ad_{Y}^{\ell-j}\; \ad_{X}^{j}Y\| \nonumber\\
 &\leq  &  \frac{1}{\ell!\;2^{\ell+1}} 2^{\ell}\, y^{\ell}\, x+ \frac{1}{2^{\ell}}\sum_{j=0}^{\ell}\,\frac{1}{j!(\ell-j)!} y^{\ell-j} \, x^j\,y \, 2^{\ell-j}\,2^j  \nonumber\\
 & = & \frac{1}{2\,\ell!} \, y^{\ell}\, x+ \sum_{j=0}^{\ell} \frac{1}{j!\, (\ell-j)!} \,y^{\ell+1-j}\, x^j.
\end{eqnarray}
Since
$$
(x+y)^k=\sum_{j=0}^{k} \, {k \choose j} \, x^j \, y^{k-j}\; = \; \sum_{j=0}^{k}\, \frac{k!}{j!\,(k-j)!} \, x^j \, y^{k-j},
$$
then
\begin{eqnarray}    \label{cota2:f1k}
 \|f_{1,\ell}\| & \leq & \frac{1}{\ell!}\; \left( \frac 12 \, y^{\ell}\, x+y\, (x+y)^{\ell}\right) \\
  &  \le &  \frac{1}{\ell!} (x + y)^{\ell+1} \equiv s_{1,\ell} (x+y)^{\ell+1}.  \nonumber
 \end{eqnarray}
Analogously, from (\ref{eq:G_1.1}) it is clear that $\|g_{1,0}\| \le (x+y)/2$ and
\begin{equation}   \label{bound.g}
  \|g_{1,\ell}\| \le \frac{1}{2 \ell!} y^{\ell} x \le \frac{1}{2 \ell!}  (x+y)^{\ell+1} \equiv \tilde{s}_{1,\ell} (x+y)^{\ell+1},
\end{equation}
so that
\[
  \|C_3\| \le \frac{1}{6} (\|f_{1,2}\| + \|g_{1,2}\|) \le \frac{1}{8} (x+y)^3 \equiv r_3 (x+y)^3.
\]
Proceeding by induction according with algorithm (\ref{alg.1}) one has in general
\[
   \|f_{k,\ell}\| \le s_{k,\ell} (x+y)^{\ell+1}, \qquad \|g_{k,\ell}\| \le \tilde{s}_{k,\ell} (x+y)^{\ell+1}, \qquad \|C_k\| \le r_k (x+y)^k
\]
with
\begin{equation}  \label{bounds.fg}
\aligned
 & s_{1,\ell}=\frac{1}{\ell!}, \qquad
  s_{k,\ell}=\left\{\begin{array}{ll}
                      \displaystyle  \sum_{j=0}^{[\ell/k]}  \frac{1}{j!} 2^j \, r_{k}^j s_{k-2,\ell-kj},& \  \ell \neq k-1 \\[5mm]
                        k \; r_{k}+ s_{k-2,k-1}, & \   \ell=k-1
                      \end{array} \right., \\
 & \tilde{s}_{1,\ell}=\frac{1}{2 \ell!}, \qquad
  \tilde{s}_{k,\ell}=\left\{\begin{array}{ll}
                      \displaystyle  \sum_{j=0}^{[\ell/k]}  \frac{1}{j!} 2^j \, r_{k}^j \tilde{s}_{k-2,\ell-kj},& \  \ell \neq k-1 \\[5mm]
                        k \; r_{k}+ \tilde{s}_{k-2,k-1}, & \   \ell=k-1
                      \end{array} \right. \\
  & r_k = \frac{1}{2k} (s_{k-2,k-1} + \tilde{s}_{k-2,k-1})
\endaligned
\end{equation}
The series (\ref{conv.exp}) converges if the power series $\sum_{k \ge 3} \lambda^k r_k (x+y)^k$ does, and a sufficient condition is given by
\[
  \lim_{k \rightarrow \infty} \frac{r_{k+2} \lambda^{k+2} (x+y)^{k+2}}{r_k \lambda^k (x+y)^k}  =  \lambda^2 (x+y)^2 \lim_{k \rightarrow \infty} \frac{r_{k+2}}{r_k}  < 1.
\]
Convergence is then ensured as long as
 \begin{equation}  \label{eq.esti.1}
   \lambda (x+y) < \frac{1}{\sqrt{r}},
\end{equation}
where $r \equiv  \lim_{k \rightarrow \infty}  \frac{r_{k+2}}{r_k}$. Computing the coefficients $r_k$ according with the recursion (\ref{bounds.fg}) for
sufficiently large values of $k$ gives $r \approx 0.5717$ and so the symmetric Zassenhaus formula converges if $ \lambda (x+y) < 1.3225$.

\begin{figure}[ht]
\begin{center}
  \includegraphics[scale=1.3]{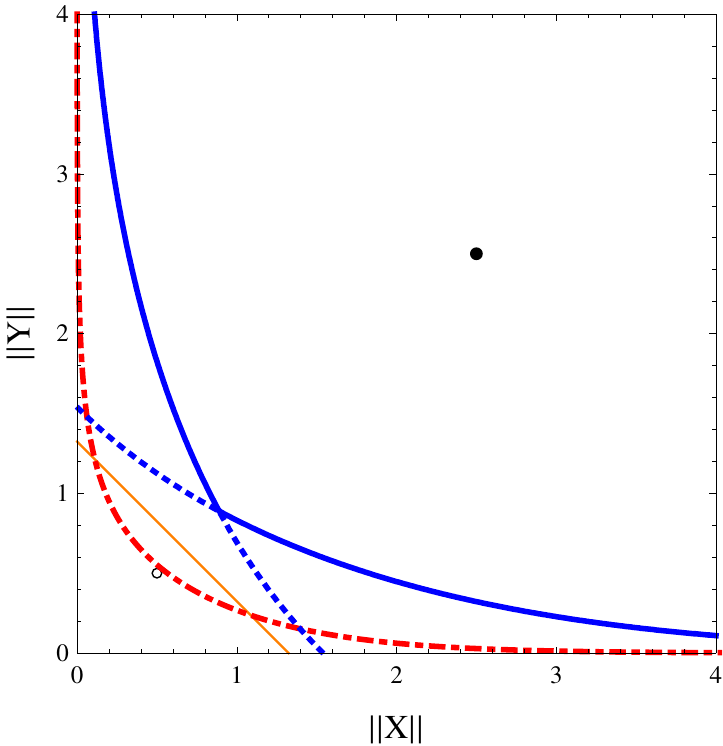}
\end{center}
\caption{{\small Bounds obtained for the convergence domain of the symmetric Zassenhaus formula: $\mathcal{D}_1 \cup \mathcal{D}_2$ (thick solid line)
and region $\|X\| + \|Y\| < 1.3225$ (thin solid line). The result obtained in  \cite{casas12eco} for the standard Zassenhaus formula (\ref{zass.1}) is also
included for comparison (dash-dotted curve).}}
\label{fig:1}
\end{figure}

A better
estimate can be obtained, however, by taking sharper bounds for $\|f_{1,\ell}\|$ and $\|g_{1,\ell}\|$ and then applying again algorithm (\ref{alg.1}).
Specifically, from (\ref{cota2:f1k}) and (\ref{bound.g}) we get, respectively
\[
   \|f_{1,\ell}\| \le d_{1,\ell}, \qquad\quad \|g_{1,\ell}\| \le \tilde{d}_{1,\ell},
\]
with $d_{1,0} = \tilde{d}_{1,0} = \frac{1}{2} (x+y)$ and
\begin{equation}   \label{d1ele}
   d_{1,\ell} = \frac{1}{\ell!}\; \left( \frac 12 \, y^{\ell}\, x+y\, (x+y)^{\ell}\right), \qquad\quad
  \tilde{d}_{1,\ell} = \frac{1}{2 \ell!} y^{\ell} x, \qquad \ell \ne 0.
\end{equation}
Assume for clarity that $\lambda=1$. Then, from the recursion (\ref{alg.1}) we have for $k=3,5,\ldots$
\begin{eqnarray}    \label{cota:fnk}
 \|f_{k,\ell}\| & \leq & d_{k,\ell} \equiv \sum_{j=0}^{[\ell/k]}\frac{1}{j!}\, 2^j \, \delta_k^j \, d_{k-2,\ell-kj}, \qquad d_{k0}=d_{10} \nonumber\\
  \|f_{k,k-1}\| & \leq & d_{k,k-1} \equiv k\, \delta_k\; + \; d_{k-2,k-1} \\
 \|g_{k,\ell}\| & \leq & \tilde{d}_{k,\ell} \equiv \sum_{j=0}^{[\frac \ell/k]}\frac{1}{j!}\, 2^j \, \delta_k^j \, \tilde{d}_{k-2,\ell-kj}, \qquad \tilde{d}_{k0}=\tilde{d}_{10} \nonumber\\
  \|g_{k,k-1}\| & \leq & \tilde{d}_{k,k-1} \equiv k\, \delta_k\; + \; \tilde{d}_{k-2,k-1} \nonumber  \\
 \|C_k\| & \leq & \delta_k \equiv \frac{1}{2k}\left(d_{k-2,k-1} +\tilde{d}_{k-2,k-1}\right).  \nonumber
\end{eqnarray}
As before, if the series $\sum_{k \ge 3} \delta_k$ converges, so does $M(\lambda=1)$. Therefore,
 sufficient condition for convergence of the symmetric Zassenhaus formula is obtained by imposing
\begin{equation}   \label{conv}
   \lim_{k \rightarrow \infty}  \frac{\delta_{k+2}}{\delta_k} < 1.
\end{equation}
At this point it is worth remarking that, although not reflected by the notation,
$d_{k,\ell}$, $\tilde{d}_{k,\ell}$ and $\delta_k$ all depend on $(x,y)=(\|X\|,\|Y\|)$, and so condition (\ref{conv}) implies in fact a
constraint on the convergence domain $(x,y) \in \mathbb{R}^2$ of the symmetric Zassenhaus formula.
In Figure~\ref{fig:1}, we depict the (numerically computed) domain $\mathcal{D}_1$ of such points $(x,y)$.
This has been obtained by considering a sufficiently fine partition in the $x$-axis and for each value of $x$
we have determined the highest value of $y$ verifying condition (\ref{conv}) by computing the coefficients $d_{k,\ell}$ and $\delta_k$ up to $k=401$
(although considering a smaller value of $k$ the graph does no change significantly).
 The domain of convergence contains the point (0.001, 1.539)
 and also the points $(x,0)$ with arbitrary large values of $x$. On the other
hand, since the roles of $X$ and $Y$ are interchangeable (equation (\ref{zassym.red2})), we can obtain similar bounds for the terms $\tilde{C}_i$
simply by applying the transformation $x \leftrightarrow y$ in (\ref{d1ele}) and (\ref{cota:fnk}), thus resulting in the domain $\mathcal{D}_2$ (symmetric
of $\mathcal{D}_1$). The boundary of $\mathcal{D}_1 \cup \mathcal{D}_2$ corresponds to the thick solid line in Figure~\ref{fig:1}.
For completeness we have also included the estimate
(\ref{eq.esti.1}), i.e., the region $x+y < 1.3225$ (thin solid line) and also the convergence domain for the standard Zassenhaus formula (\ref{zass.1})
obtained in \cite{casas12eco} (dash-dotted curve). Notice the great improvement obtained when considering the symmetric version of the formula.

\section{Examples}

In this section we illustrate the behavior of the symmetric Zassenhaus formula in comparison with the
standard Zassenhaus formula on several examples of physical and mathematical interest.

\paragraph{Example 1.} In reference \cite{goldman14pdq} several approximations for the effective Hamiltonian in periodically driven quantum systems are constructed. In particular,
the Zassenhaus formula is used to find factorized expressions for the individual operators involved in a particular splitting of the evolution operator when the different
pieces of the Hamiltonian satisfy the following cyclic relations:
\begin{equation}  \label{ex.1.1}
  [X,Y] = Z, \qquad [X,Z] = a Y, \qquad [Y,Z]=0,
\end{equation}
where $a$ is a parameter. Then, a direct application of (\ref{zass.1}) leads to \cite{goldman14pdq}
\begin{equation}   \label{ex.1.2}
   \e^{X + Y} = \e^X \, \e^{f_1(a) Y + f_2(a) Z}
\end{equation}
with
\[
   f_1(a) = \frac{1}{\sqrt{a}} \sinh \sqrt{a}, \qquad f_2(a) = \frac{1 - \cosh \sqrt{a}}{a}.
\]
We can apply of course algorithm (\ref{alg.1}) and determine the successive terms $C_k$. Then we find
\[
   C_3 = \frac{1}{48} a Y, \qquad C_5 = \frac{1}{3840} a^2 Y, \qquad C_7 = \frac{1}{645120} a^3 Y, \ldots
\]
so that   
\[
   \e^{X + Y} = \e^{\frac{X}{2}} \, \exp \left( Y + \frac{1}{24} a Y + \frac{1}{1920} a^2 Y + \frac{1}{322560} a^3 Y + \cdots \right) \, \e^{\frac{X}{2}}
\]
and we arrive at the closed form result
\begin{equation}   \label{ex.1.3}
       \e^{X + Y} = \e^{\frac{X}{2}} \, \e^{f_3(a) Y} \, \e^{\frac{X}{2}}, \qquad \mbox{ with } \qquad f_3(a) = \frac{2}{\sqrt{a}} \sinh \left( \frac{\sqrt{a}}{2} \right).
\end{equation}
It is worth noticing that a simpler expression for $\exp(X+Y)$ is obtained in this case which is also independent of the $Z$ operator.

\paragraph{Example 2.}  A Lie algebra of interest for quantum mechanical problems related with the harmonic oscillator is spanned by the operators
$\{Q, P, W \equiv P^2 + Q^2, cI\}$ with the commutation relations  \cite{wilcox67eoa}
\[
      [P, Q] = c I, \quad [W,P] = -2 c Q, \quad [W,Q] = 2 c P, \quad [P,I] = [Q,I]= [W, I]=0.
\]
An alternative set of operators is formed by $\{W,X,Y, sI \}$, where 
\[
  X = Q -i P, \quad Y = Q + i P, \quad W= X Y - \frac{1}{2} s I, \quad s = 2 i c,
\]
verifying
\begin{equation}  \label{ex.2.1}
   [W, X] = -s X, \quad [W, Y] = s Y, \quad [X, Y] = s I.
\end{equation}
In this case, the symmetric Zassenhaus formula reduces to 
\[
    \e^{X+Y} = \e^{\frac{1}{2} X} \, \e^Y \,  \e^{\frac{1}{2} X}.         
\]
On the other hand, algorithm (\ref{alg.1}) also provides a closed expression for $\exp(X+W)$. Specifically,
\[
  \e^{X + W} = \e^{\frac{1}{2} X} \, \e^{\frac{1}{2} W} \, \e^{g(s) X} \, \e^{\frac{1}{2} W} \,  \e^{\frac{1}{2} X},
\]
where the function $g(s)$ is defined through the power series
\begin{equation}  \label{def.g}
   g(s) =   \sum_{k=1}^{\infty}  d_{2k+1} s^{2k}
\end{equation}
whose first coefficients read
\[ 
   d_3 = -\frac{1}{12}, \quad d_5 = -\frac{1}{480}, \quad d_7 = -\frac{1}{53760}, \quad d_9 = -\frac{1}{11612160}.
\]
More in general, one has
\[
   \e^{X + Y + W} =   \e^{\frac{1}{4} X}  \,  \e^{\frac{1}{2} Y} \,  \e^{\frac{1}{4} X} \,  \e^{\frac{1}{2} W} \, V \,  \e^{\frac{1}{2} W} \,    \e^{\frac{1}{4} X} \, 
     \e^{\frac{1}{2} Y} \, \e^{\frac{1}{4} X}, 
\]
with
\[
   V = \e^{-h(s)} \, \e^{\frac{1}{2} g(s) X} \, \e^{g(s) Y} \, \e^{\frac{1}{2} g(s) X}.
\]
Here $g(s)$ is the function introduced in (\ref{def.g}) and 
\begin{equation}  \label{def.h}
   h(s) =   \sum_{k=1}^{\infty}  d_{2k+1} 4^{k-1} (2k-3) s^{2k}.
\end{equation}
Other closed form expressions for $\e^{X+Y+W}$ are obtained, of course, by interchanging the role of the different operators in algorithm (\ref{alg.1}).

\begin{figure}[h!]
\begin{center}
  \includegraphics[scale=1]{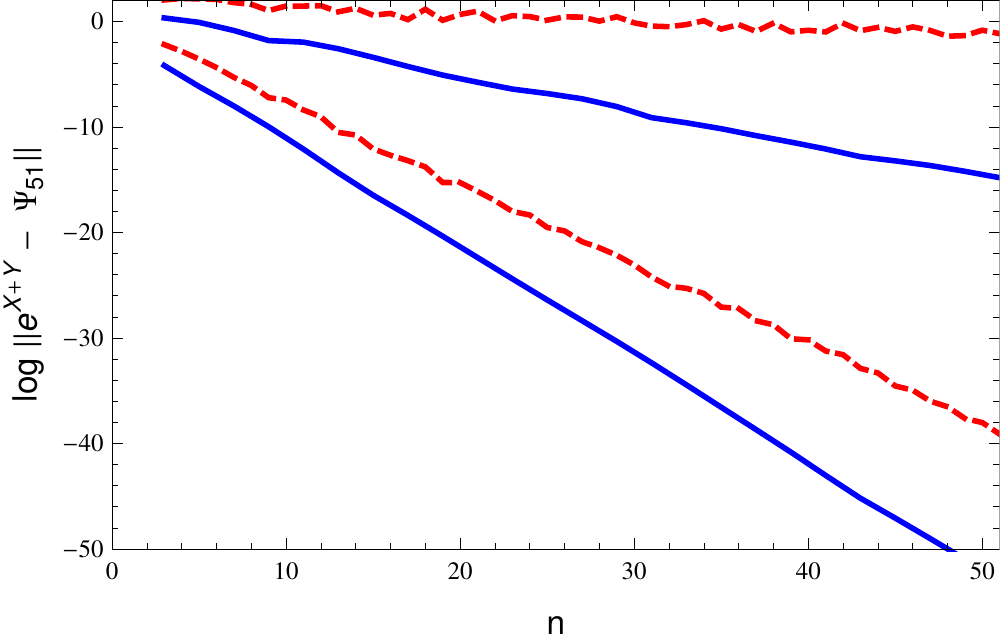}
\end{center}
\caption{\small{Error (in logarithmic scale) of the approximation
of order $n=51$ for random $20 \times 20$ matrices with $\|X\| =\|Y\|=2.5$ (top curves) and with $\|X\| =\|Y\|=0.5$ (bottom) obtained with the standard Zassenhaus
formula (dashed lines) and its symmetric version (solid lines).}}
\label{fig:2}
\end{figure}

\paragraph{Example 3.}  Our next example illustrates in practice the differences between the usual Zassenhaus formula and its symmetric
version, both in the domain of convergence and in their rate of convergence. To do that, we take 
two pairs of random  $20 \times 20$ matrices $X$ and $Y$. The first couple is chosen in such
a way that the norm of $X$ and $Y$ both belong to the subset of $\mathcal{D}_1 \cup \mathcal{D}_2$ within the convergence domain of the
standard Zassenhaus formula. Thus, we take $x=y=0.5$, compute the first $51$ terms of both factorizations (\ref{zassym.red}) and (\ref{zass.1})
 and finally determine the norm of the error of the respective approximations. More specifically, we compute $\|\e^{X+Y} - \Psi_{n}\|$, where
\begin{equation}  \label{psi.sym}
   \Psi_{n} = \e^{\frac{1}{2} X} \, \e^{\frac{1}{2}  Y} \,
   \e^{ C_2} \, \e^{ C_3} \,  \cdots   \e^{2 C_{n}}  \,  \cdots \e^{ C_3} \,
   \e^{C_2} \, \e^{\frac{1}{2}  Y} \, \e^{\frac{1}{2} X}
\end{equation}
with the $C_i$ generated by (\ref{alg.1}) and
\begin{equation}  \label{psi.ns}
   \Psi_{n} = \e^{X} \, \e^{Y} \,
   \e^{C_2} \, \e^{C_3} \,  \cdots  \e^{C_{n}},
\end{equation}
with the terms $C_i$ obtained in \cite{casas12eco}. In Figure \ref{fig:2} we depict this error as a function of the number of terms $n$ in the approximation
(\ref{psi.sym}) corresponding to the symmetric Zassenhaus formula (solid curve) and (\ref{psi.ns}) for the standard non-symmetric version (dashed line).
Although the convergence of both procedures is clearly visible from the graphs in this case,
the rate of convergence of the symmetric approximation is clearly faster.

For the second pair we take matrices $X$ and $Y$ such that $x =y=2.5$, and repeat the experiment, thus providing the top curves in Figure \ref{fig:2}. Although
we are now clearly outside the guaranteed convergence domain of the symmetric Zassenhaus formula,
we observe that the approximations (\ref{psi.sym}) still converges to the exact result (although at a smaller rate of convergence). This is not the case of
the standard Zassenhaus formula, where the error does not diminish significantly with $n$. We see then that our estimate for the convergence
of the symmetric Zassenhaus formula is not sharp: there are still matrices with larger norm for which one gets convergent approximations.

\

For convenience, we have also depicted in Figure \ref{fig:1} the exact location of both pair of matrices by an empty and a full circle, respectively.

\begin{figure}[h!]
\begin{center}
  \includegraphics[scale=1]{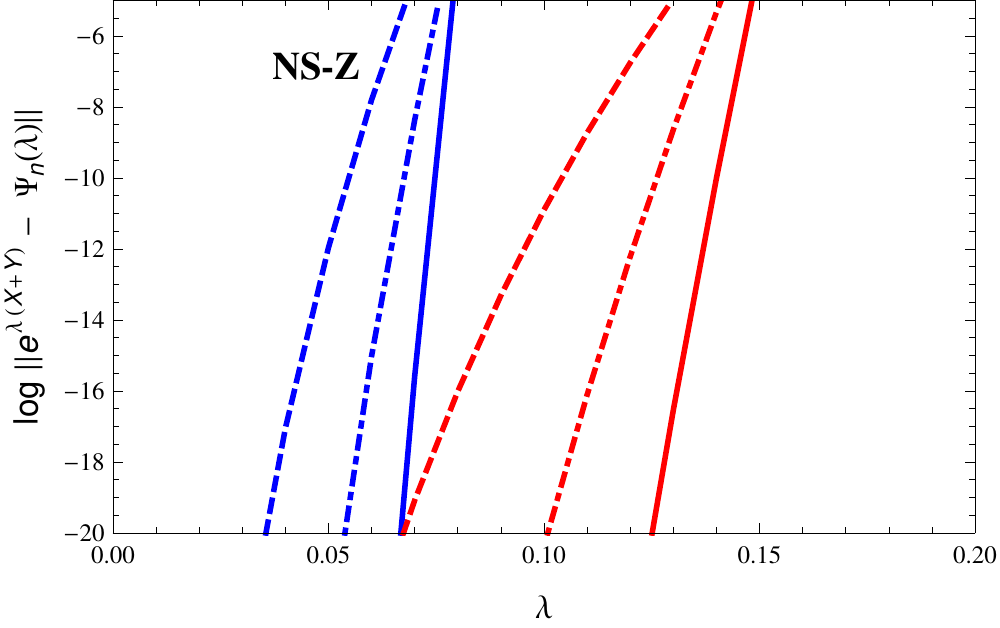}
\end{center}
\caption{\small{Error in the approximation $\Psi_n(\lambda)$ corresponding to (\ref{psi.sym}) and (\ref{psi.ns}) for several values of $n$: 51 (dashed), 101 (dash-dotted) and 201 (solid).
Curves labelled by NS-Z correspond to the standard Zassenhaus formula.}}
\label{fig:3}
\end{figure}

\paragraph{Example 4.} For our final illustration we take $X$ and $Y$ to be  $2 \times 2$ matrices and analyze the approximation $\Psi_n(\lambda)$ in
(\ref{conv.2}) obtained with the recurrence (\ref{alg.1}) as a function of $\lambda$ for a given number of terms $n$. Specifically, we consider
\[
  X = \pi \left( \begin{array}{cr}
    		0  & \  \alpha \\
		-1/\alpha  &  0
	  \end{array}  \right), \qquad\quad
  Y = \pi \left( \begin{array}{cc}
    		0  & \  (10 + 4 \sqrt{6}) \alpha \\
		(-10+4 \sqrt{6})/\alpha  &  0
	  \end{array}  \right).	  	
\]
It turns out that $\e^{X+Y} = \e^X \, \e^Y = \e^Y \, \e^X$ for all $\alpha \ne 0$, although $[X,Y] \ne 0$ \cite{frechet52lsn}, but $\e^{\lambda(X+Y)} \ne \e^{\lambda X} \, \e^{\lambda Y}$
for $\lambda \ne 1$. We take $\alpha=1/5$, so that
$x = 15.7205 $, $y=12.8379$ (in the 2-norm), and depict the error $\|\e^{\lambda(X+Y)} - \Psi_{n}(\lambda)\|$ for several values of $n$ in the approximation (\ref{conv.2})
vs. the parameter $\lambda$. Specifically, in Figure \ref{fig:3} we collect the results achieved when $n=51$ (dashed line), $n=101$ (dash-dotted line) and
$n=201$ (solid line). We also carry out the analogous computation for the standard Zassenhaus formula (curves labelled by NS-Z in the figure). Notice the remarkable improvement
obtained when considering the symmetric approximation. Here, once again, we observe a convergent behavior for values of $x$ and $y$ well beyond the
domain $\mathcal{D}_1 \cup \mathcal{D}_2$ in Figure \ref{fig:1}. In particular, for $\lambda=0.13$ we still have a valid approximation within 18 digits when $n=201$.

\

In summary, we have analyzed a symmetric version of the Zassenhaus formula for disentangling the exponential of the sum of two
noncommuting operators, eq. (\ref{zass.1.2}). We have developed a recursive procedure that allows one in principle to get all the terms
in the factorization up to a prescribed degree.
Compared with the usual Zassenhaus formula, eq. (\ref{zass.1}), it possesses a larger domain of convergence in the matrix case and its rate
of convergence is also faster. Moreover, the examples considered show that, when the whole series can be computed in closed form, the new procedure is able to provide more compact 
expressions, whereas in the general case, if the factorization is truncated at a certain degree $n$ then one gets better approximations than with (\ref{zass.1}) and at the same time the number of
terms to compute is smaller (all terms for even $n$ vanish). In consequence, one can advantageously replace (\ref{zass.1}) by its symmetric
version (\ref{zass.1.2}) in real applications.

\subsection*{Acknowledgments}
The authors would like to thank the anonymous referee for his/her comments that have contributed to clarify some points in the presentation in a previous version
of the manuscript.
This work has been funded by the research project P1.1B20115-16 (Universitat Jaume I). FC has been additionally supported
by Ministerio de Econom\'{\i}a y Competitividad (Spain) through  project MTM2016-77660-P (AEI/FEDER, UE).




\end{document}